\documentclass[11pt]{article}
\usepackage{setspace}
\usepackage{amssymb,latexsym,amsmath,amsbsy,amsthm,amsxtra,amsgen}

\doublespace
\newfont{\msbm}{msbm10 at 11pt}

\def\mn{\bigskip\noindent}

\newtheorem{Theo}{Theorem}
\newtheorem{Lemma}[Theo]{Lemma}

\begin{document}
\title{A stochastic model for phylogenetic trees}
\author{by Thomas M. Liggett\thanks{Partially supported by NSF grant DMS-0301795} and Rinaldo B. Schinazi\thanks{Partially supported by NSF grant DMS-0701396}\\\\
University of California at Los Angeles,\\
and University of Colorado at Colorado Springs \\
}
\maketitle

\begin{abstract}
We propose the following simple stochastic model for phylogenetic trees. 
New types are born and die according to a birth and death chain. At each birth we associate a fitness 
to the new type sampled
from a fixed distribution. At each death the type with the smallest fitness is killed. We show that if the
birth (i.e. mutation) rate is subcritical we get a phylogenetic tree consistent with an influenza tree (few
types at any given time and one dominating type lasting a long time). When the birth rate is supercritical we
get a phylogenetic tree consistent with an HIV tree (many types at any given time, none lasting very long). 
\end{abstract}

\section{Introduction}

The influenza phylogenetic tree is peculiar in that it is very skinny:
one type dominates for a long time and any other type that arises quickly
dies out. Then the dominating type suddenly dies out and is immediately replaced 
by a new dominating type. The models proposed so far are very complex and make
many assumptions. See for instance Koelle et al. (2006) and van Nimwegen (2006). 
We would like to use a simple stochastic model for such a tree. The other motivation for this work comes
from the comparison between influenza and HIV phylogenetic trees. An HIV tree is characterized by a radial
spread outward from an ancestral node, in sharp contrast with an influenza tree. Moreover, Korber et al.
(2001) note that the influenza virus is less diverse worldwide than the HIV virus is in Amsterdam alone.
However, both types of trees are supposed to be produced by the same basic mechanism: mutations. Can the same
mathematical model produce two trees that are so different? Our simple stochastic model will show a striking
difference in behavior depending on the mutation rate.

Our model has a birth and death component and a fitness component. For the death and birth component we do the following. If there are $n\geq 1$ types at a certain time $t$ then there is birth of a new type (by mutation) at rate $n\lambda$.
We think of a birth as the appearance of one new type, not the replacement of one type by two new types. If there are $n\geq 2$ types then there is death of one type at rate $n$. If only one type is left it cannot die. 
That is,
\begin{align}
n&\longrightarrow n+1\text { at rate }n\lambda
\nonumber \\
n&\longrightarrow n-1\text { at rate }n\text{ if }n\geq 2. \nonumber
\end{align}

Moreover, each new individual is assigned a fitness value chosen from a fixed distribution, independently each time. Every time there is a death event then the
type that is killed is the one with the smallest fitness. Since all that matters is the ranks of the fitnesses, we might as well take their distribution to be uniform on $[0,1]$. For simplicity the process is started with a single type.

We give no specific rule
on how to attach a new type after a birth to existing types (in order to construct a tree). Our results do not
depend on such a rule. Two natural possibilities are to either attach the new type to the type which has the
maximum fitness or to a type taken at random. 
\begin{Theo} 
 Take $\alpha\in (0,1)$. 

\noindent If $\lambda\leq 1$, then 
$$\lim_{t\rightarrow\infty}P(\text{maximal types at times } \alpha t \text{ and }t\text{ are the same})
=\alpha,$$ while if $\lambda>1$, then this limit is 0.
\end{Theo}
We see that if $\lambda< 1$, the dominating type (i.e. the fittest type) at time $t$ has likely been present
for a time of order $t$ and at any given time there will not be many types. This is consistent with 
the observed structure of an
influenza tree. On the other hand, if $\lambda>1$, then the dominating type at time $t$ has likely been present
for a time of order shorter than $t$ and at any given time there will be many types. This is consistent with
an HIV tree. 

\section{Proof of Theorem 1}

The proof divides into three cases, depending on whether the birth and death chain is positive
recurrent, null recurrent, or transient. We present them in order of difficulty.
\subsection{Case $\lambda<1$}
Let $\tau_1, \tau_2,...$ be the (continuous) times
between successive visits of the chain to 1, $T_n=\tau_1+\cdots + \tau_n$, $\sigma_1,\sigma_2,...$ be the
number of new types introduced in cycles between successive visits to 1, and
$S_n=1+\sigma_1+\cdots+\sigma_n$. Note that
the $\tau$'s and
$\sigma$'s are not independent of each other, but the sequence $ (\tau_1,\sigma_1),(\tau_2,\sigma_2),...$ is i.i.d.
and independent of the fitness sequence. Define the usual renewal process $N(t)$ corresponding to
the $\tau$'s by $\{N(t)=n\}=\{T_n\leq t<T_{n+1}\}$.

For $0<s<t$, recalling that $T_{N(t)}\leq t< T_{N(t)+1}$, and noting that the
maximal type is increasing in time, we see that

\begin{equation}
P(\text{maximal types at times }s \text{ and }t\text{ are the same},N(s)<N(t))
\label{(1)}
\end{equation}

\noindent  lies between
$$P(\text{maximal types at times }T_{N(s)} \text{ and }T_{N(t)+1}\text { are the same},N(s)<N(t))$$
and
$$P(\text{maximal types at times }T_{N(s)+1} \text{ and }T_{N(t)}\text { are the same},N(s)<N(t)).$$
Let ${\cal F}$ be the $\sigma$-algebra generated by
$(\tau_1,\sigma_1),(\tau_2,\sigma_2),...$. Then for $k\leq l$, since the 
fitness sequence is i.i.d. and independent of ${\cal F}$,
$$P(\text{maximal types at times }T_k \text{ and }T_l\text{ are the same} \mid  {\cal F})=\frac{S_k}{S_l}.$$
More precisely, conditional on ${\cal F}$ there are $S_l$ fitnesses observed by time $T_l$ and $S_k$ of them are observed by time $T_k$. We claim that in $n$ i.i.d. observations the probability that the largest occurs among the first $m$ is $m/n$, since any one of $n$ is equally likely to be the largest. Since $N(s)$ and $N(t)$ are $ {\cal F}$ measurable, it follows that (1) lies between
\begin{equation}
E\bigg[\frac{S_{N(s)}}{S_{N(t)+1}},N(s)<N(t)\bigg]\quad\text{and}\quad
E\bigg[\frac{S_{N(s)+1}}{S_{N(t)}},N(s)<N(t)\bigg].
\label{(2)}
\end{equation}

Since $\lambda<1$, $E\tau<\infty$, and the renewal theorem gives 
$$N(s)/s\rightarrow 1/E\tau\text{ a.s.,}$$ while the strong law of large numbers gives
$S_{N(s)}/N(s)\rightarrow E\sigma$ a.s., so that $S_{N(s)}/s \rightarrow E\sigma/E\tau$ a.s.
It follows by the bounded convergence theorem that
\begin{equation}
\lim_{t\rightarrow\infty}P(\text{maximal types at times } \alpha t \text{ and }t\text{ are the same})
=\alpha.
\label{(3)}
\end{equation}
This completes the proof of Theorem 1 in the subcritical case.
\subsection{Case $\lambda> 1$.}
Define the $\tau$'s and $\sigma$'s as above, except that now, the cycles used are between the successive
times the chain reaches a new high. In other words, $T_n$ is the hitting time of $n+1$, 
$\sigma_n$ is the number of new types born during a first passage cycle from $n$ to $n+1$ and 
$S_n$ is the number of new types seen up to time $T_n$. Of course, the $\sigma$'s and $\tau's$
are no longer identically distributed. However, $(\tau_1,\sigma_1), (\tau_2,\sigma_2),...$ are independent. 
The key to the proof is the following Lemma.
\begin{Lemma}
Assume that $\lambda>1$. Then $e^{-(\lambda-1)t}N(t)$ is almost surely bounded. 
\end{Lemma}
\begin{proof}[Proof of Lemma 2]
Our first step in this proof is to estimate the first two moments of $\tau_n$. Following Keilson (1979) (see (5.1.2)) we note that $\tau_n$ has the same distribution as 
\begin{equation}
{X\over (1+\lambda)n}+Y(\tau_{n-1}+\tau_n')\text{ for }n\geq 2,
\label{(4)}
\end{equation}
where $X$ has a mean 1 exponential distribution, $\tau_n'$ has the same distribution as $\tau_{n}$, $Y$ is a Bernoulli random with $P(Y=1)={1\over \lambda+1}$, and $X$, $Y$, $\tau_{n-1}$ and $\tau_n'$ are independent. Letting $\mu_n=E\tau_n$, it follows from (4) that
\begin{equation}
\lambda\mu_n={1\over n}+\mu_{n-1}\text{ for $n\geq 2$ and }\mu_1={1\over\lambda}.
\label{(5)}
\end{equation}
We will use the following recursion formula, which is easy to prove by induction.
\begin{Lemma}
Let $a_n$ and $b_n$ be two sequences of real numbers such that
$a_1=\lambda^{-1}b_1$ and  for $n\geq 2$, $\lambda a_n=b_n+a_{n-1}$. Then,
$$a_n=\sum_{j=1}^n \lambda ^{-j}b_{n+1-j}.$$
\end{Lemma}

Applying Lemma 3 to (5) we get
$$\mu_n=\sum_{j=1}^n {\lambda ^{-j}\over n+1-j}.$$
Writing $1/(\lambda-1)$ as a geometric series we have
$$\mu_n-{1\over n(\lambda-1)}=\sum_{j=1}^n \lambda ^{-j}({1\over n+1-j}-{1\over n})
-{1\over n}\sum_{j=n+1}^\infty \lambda ^{-j}.$$
Changing the order of summation gives
$$\sum_{n=1}^\infty \sum_{j=1}^n \lambda ^{-j}({1\over n+1-j}-{1\over n})=\sum_{j=1}^\infty\lambda ^{-j}\sum_{n=j}^\infty ({1\over n+1-j}-{1\over n}).$$
Note that for $j\geq 2$
$$\sum_{n=j}^\infty ({1\over n+1-j}-{1\over n})=\sum_{k=1}^{j-1}{1\over k}$$
and this term is 0 for $j=1$.
Hence,
$$\sum_{n=1}^\infty {1\over n}\sum_{j=n+1}^\infty \lambda ^{-j}=\sum_{j=2}^\infty\lambda ^{-j}
\sum_{k=1}^{j-1}{1\over k}.$$
We conclude that
$$\sum_{n=1}^\infty|\mu_n-{1\over n(\lambda-1)}|<\infty$$
and
$$\sum_{n=1}^\infty(\mu_n-{1\over n(\lambda-1)})=0.$$
Therefore,
\begin{equation}
E(T_n)-{1\over\lambda -1}\sum_{k=1}^n {1\over k}\text{ converges  to }0.
\label {(6)}
\end{equation}

We also need an almost sure result for $T_n$, and for this, we will estimate the second moment
of $\tau_n$. Let $v_n=Var(\tau_n)$. It is easy to check that if
$Y$ is a Bernoulli random variable and is independent of a random variable $Z$ then
$$Var(ZY)=E(Y)Var(Z)+Var(Y)(EZ)^2.$$
Using this remark and (4) we have for $n\geq 2$
$$
v_n={1\over (1+\lambda)^2n^2}+{1\over 1+\lambda}(v_n+v_{n-1})+{\lambda\over(1+\lambda)^2}(\mu_n+\mu_{n-1})^2
$$
Therefore, for $n\geq 2$
\begin{equation}
\lambda v_n=b_n+v_{n-1}
\label {(7)}
\end{equation}
where
$$b_n={1\over (1+\lambda)n^2}+{\lambda\over 1+\lambda}(\mu_n+\mu_{n-1})^2.$$
Set $\mu_0=0$, then $\lambda v_1=b_1$. Hence, Lemma 3 applies to (7), giving
\begin{equation}
v_n=\sum_{j=1}^n \lambda ^{-j}b_{n+1-j}.
\label {(8)}
\end{equation}
Since $\mu_n\sim {1\over n(\lambda-1)}$ (that is, the ratio converges to 1), 
$b_n\sim {C\over n^2}$ where $C$ depends on $\lambda$ only. From (8) we get
$$v_n\sim C'b_n\sim {CC'\over n^2},$$
where $C'$ depends on $\lambda$ only. 
This implies the a.s. convergence of the
random series $\sum(\tau_n-E\tau_n)$
(see for instance Corollary 47.3 in Port (1994)). Therefore the partial sums converge a.s. and
$$T_n-E(T_n)\text { converges a.s.}$$
Using (6) we get that
$T_n-{1\over \lambda -1}\log n$ converges a.s. and is therefore a.s. bounded. Now use the fact that $\{N(t)\geq n\}=\{T_n\leq t\}$ to conclude that $N(t)\exp(-(\lambda-1)t)$ is almost surely bounded. This concludes the proof of Lemma 2.
\end{proof}
We are now ready to complete the proof of Theorem 1 in the supercritical case. 
Let $(Z_{i})_{i\geq 1}$ be a discrete time random walk starting at 0 that goes to the right with probability $\lambda/(\lambda+1)$ and
to the left with probability $1/(\lambda+1)$. For every $n\geq 1$, let $Z_{i,n}$ be
a discrete time random walk starting at 0 with the same rules of evolution as $Z_i$ except that the random walk $Z_{i,n}$ has a reflecting barrier at $-n+1$. For every $n\geq 1$, the two random walks $Z_i$ and $Z_{i,n}$ are coupled so that they move together until (if ever) they hit $-n+1$ and thereafter we still couple them so that 
$Z_i\leq Z_{i,n}$ for every $i\geq 0$.
Let $U$ and $U_n$ be the hitting times of 1 for the random walks $Z_i$ and $Z_{i,n}$, respectively. 

First note that a new type appears every time there is a birth. Therefore, $\sigma_n$ is the number of steps to the right of the random walk $Z_{i,n}$ stopped at 1. That is, $\sigma_n$ is $(1+U_n)/2$. We now show that $U_n$ converges a.s to $U$. 
Let $\delta>0$ we have
$$P(|U_n-U|>\delta)\leq P(U>U_n)\leq P(Z_i=-n+1 \text{ for some }i\geq 1).$$
The last probability decays exponentially with $n$. Therefore,
$$\sum_{n\geq 1}P(|U_n-U|>\delta)<\infty.$$
An easy application of Borel-Cantelli Lemma implies that $U_n$ converges a.s. to $U$.
Since $U_n\leq U$ the Dominated Convergence Theorem implies that, for every $k\geq 1$ the $k$th moment of $\sigma_n$ converges to the $k$th moment of $(1+U)/2$. In particular, $Var (\sigma_n)$ is a bounded sequence. This is enough to prove that
$${1\over n}\sum_{i=1}^n(\sigma_i-E(\sigma_i))\text{ converges a.s. to 0;}$$
see for instance Proposition 47.10 in Port (1994).
Since $E(\sigma_n)$ is a convergent sequence we get that $S_n/n$ converges a.s. to the limit of $E(\sigma_n)$.

Since $N(t)\rightarrow\infty$ a.s., this strong law of large numbers gives that $S_{N(t)}/N(t)$ converges to
the limiting expectation of
$\sigma_n$. This together with Lemma 2 shows that the two terms in (2) converge to 0 when we let $s=\alpha t$
and $t$ goes to infinity. The proof of Theorem 1 in the supercritical case is complete.
\subsection{Case $\lambda= 1$.}

In this subsection we go back to the notation of section 2.1 where $T_n$ be the time of the $n$th visit of the chain to 1.

\begin{Lemma} 
Let $\lambda=1$. Then,
$$\frac {T_n}{n\log n}\rightarrow 1\text{ in probability.}
$$
\end{Lemma}
\begin{proof}[Proof of Lemma 4]

When the chain hits 1, it waits a mean 1 exponential time and then jumps to 2. Hence,
$$T_n=\sum_{i=1}^n X_i+\sum_{i=1}^n H_i$$
where the $X_i$ are independent mean 1 exponential times and $H_i$ are the hitting times of 1 starting at 2. The $H_i$ are i.i.d. with distribution function $F$. From the backward Kolmogorov equation 
$$\int_0^{F(t)}{ds\over 1+s^2-2s}=t,$$
we get
$$F(t)={t\over 1+t}.$$
We now use a weak law of large numbers, see Theorem 2 in VII.7 in Feller (1971). It is easier to redo the short proof rather than check the hypotheses of the Theorem. The key is the following consequence of Chebyshev's inequality applied to the truncated random variables:
\begin{equation}
P( |{1\over nm_n}\sum_{i=1}^n H_i-1|>\epsilon)\leq {1\over n\epsilon^2 m_n^2}s_n+n(1-F(\rho_n))
\end{equation}
where
$$m_n=\int_0^{\rho_n} tF'(t)dt\text{ and }s_n=\int_0^{\rho_n} t^2F'(t)dt,$$
see (7.13) in VII.7 in Feller (1971). We will take $\rho_n=n\sqrt{\log n}$. A little Calculus shows that
$$m_n\sim \log \rho_n\sim\log n \text{ and }s_n\sim \rho_n.$$
With our choice of $\rho_n$, $n(1-F(\rho_n))$ converges to 0 and
$${1\over nm_n}\sum_{i=1}^n H_i$$
converges to 1 in probability. This completes the proof of Lemma 4.
\end{proof}
Since the events $N(t)\geq n$ and $T_n\leq t$ are the same, it follows that
\begin{equation}
N(t)\frac{\log t}t\rightarrow 1
\label{(10)}
\end{equation}
in probability as $t\uparrow\infty$.

Now, $S_n/n^2$ converges in distribution to a one sided stable law of index $\frac 12$ (see Theorem (7.7) 
in Durrett (2004)).
By (10), it follows that $S_{N(t)}/{N(t)^2}$ also has this distributional limit. 
In fact,
$${S_{N(\alpha t)}\over N^2(t)}\text{ converges to }Y_\alpha$$
in the sense of convergence of finite dimensional distributions, where $Y_\alpha$ is a stable subordinator (increasing stable process) of index 1/2.
(Note that
independence between the $\sigma$'s and $\tau$'s is not required here, which is good since they
are highly dependent. All that is needed is that the limit in (10) is constant and that
both $S_n$ and $N(t)$ are monotone.) So, the limit in (3) is 
$$\lim_{t\to\infty}E({S_{N(\alpha t)}\over S_{N(t)}})=E({Y_{\alpha}\over Y_1})=\alpha.$$  
To check the final equality, it is enough by monotonicity to verify it for rational $\alpha$.
If $\alpha=m/n$, this boils down to the simple fact that if $V_1,...,V_n$ are
i.i.d. and positive, then
$$E\frac {V_i}{V_1+\cdots+V_n}=\frac 1n.$$
\bigskip
{\bf Acknowledgements.} We thank an anonymous referee whose remarks helped improve
the presentation of the paper and simplify several of our proofs. 
\bigskip
\bigskip
\begin{center}
{\bf {\Large References}}
\end{center}
\mn R. Durrett (2004) Probability: Theory and Examples (3rd edition). Duxbury press.

\mn W. Feller (1971) An Introduction to Probability Theory and its Applications, Volume 2 (second edition). Wiley.

\mn J. Keilson (1979) Markov Chain Models-Rarity and exponentiality. Springer-Verlag.

\mn B. Korber, B. Gaschen, K. Yusim et al. (2001) Evolutionary and immunological implications of contemporary HIV-1 variation. British Medical Bulletin 58, 19-42.

\mn K. Koelle, S. Cobey, B. Grenfell and M. Pascual (2006) Epochal evolution shapes the phylodynamics of 
interpandemic influenza A (H3N2) in Humans. Science vol. 314, 1898-1903.

\mn E. van Nimwegen (2006). Influenza escapes immunity along neutral networks.
Science vol. 314, 1884-1886.

\mn S. C. Port (1994). Theoretical Probability for Applications. Wiley.
\end{document}